\documentclass[10pt]{article}
\usepackage[utf8]{inputenc}
\usepackage{latexsym,amssymb,amsmath,url,authblk}

\def\N{\mathbb{N}}
\def\Q{\mathbb{Q}}
\def\Z{\mathbb{Z}}
\def\R{\mathbb{R}}

\def\proof{\par\noindent{\em Proof. }}
\def\eproof{\hfill{$\Box$}\bigskip}

\def\ds{\dots}
\def\sus{\subset}
\def\al{\alpha}

\def\cc{\colon}

\newtheorem{thm}{Theorem}[section]
\newtheorem{prop}[thm]{Proposition}

\newtheorem{defi}[thm]{Definition}

\title{HMC real numbers in Countable Mathematical Analysis}
\author{Martin Klazar}
\date{\today}

\begin{document}

\maketitle

\medskip
\begin{abstract}
We develop a~theory of real numbers as rational Cauchy sequences, in which any two of them, 
$(a_n)$ and $(b_n)$, are equal iff $\lim\,(a_n-b_n)=0$. We need such reals in the 
Countable Mathematical Analysis (\cite{klaz_transc}) which allows to use only hereditarily at most countable (HMC) sets.
\end{abstract}

\section{Introduction}

A~set is {\em at most countable} if it is finite or countable, where the latter means that the set is in bijection
with $\omega=\{0,1,\ds\}$. A~set is {\em uncountable} if it is not at most countable. A~set 
$x$ is {\em hereditarily at most countable}, abbreviated HMC, if for every $n\in\omega$ and every chain of sets
$$
x_n\in x_{n-1}\in\ds\in x_0=x\;,
$$
the set $x_n$ is at most countable. (By the axiom of foundation every chain of sets $x_0\ni x_1\ni\ds$ is finite.) Refer to \cite{jech_ST} for set-theoretical terminology and notions.

Does one really need uncountable sets in Mathematical Analysis and in Number Theory? For example, does one need them to prove by the identities
$$
\int_0^{+\infty}x^n\mathrm{e}^{-x}\,\mathrm{dx}=n!,\ n\in\omega\;,
$$
that the Euler number $\mathrm{e}=2.71828\ds$ is transcendental? One does not, in \cite{klaz_transc} 
we carry out the proof, due to D.~Hilbert, just with HMC sets. Thus the transcendence of $\mathrm{e}$ belongs to {\em Countable Mathematical Analysis}, abbreviated CMA, 
respectively to {\em Countable Number Theory}, abbreviated CNT, where one can use only HMC sets. 

What did we do with the fact that the integrands
$$
\{(x,\,x^n\mathrm{e}^{-x})\;|\;x\in[0,\,+\infty)\}
$$
are uncountable sets? In \cite{klaz_transc} we work with their HMC restrictions to fractions
$$
\{(x,\,x^n\mathrm{e}^{-x})\;|\;x\in\Q\wedge x\ge0\}\;.
$$ 
And what did we do with the {\em min-max principle} by which every continuous 
function $f\cc[a,b]\to\R$, where $a<b$ are real numbers, attains on the interval its minimum 
and maximum? For rational restrictions to 
$$
[a,\,b]_{\Q}:=\{x\in\Q\;|\;a\le x\le b\}
$$
the principle fails as stated, there are unbounded continuous functions from
$[a,b]_{\Q}$ to $\R$. In \cite{klaz_transc} we use a~HMC variant of the principle where compactness is replaced with uniform continuity. It is based on an extension theorem.

\begin{thm}[extensions]\label{thm_extend}
Let $x$ be a~real number, $M\sus\Q$ and $f\cc M\to\R$ be a~uniformly continuous function. Then for every sequence 
$(a_n)\sus M$ with $\lim a_n=x$, the sequence
$$
(f(a_n))=(f(a_1),\,f(a_2),\,\ds)
$$
converges to a~unique real number independent of $(a_n)$ and denoted by $f(x)$. 
\end{thm}
If there exists a~sequence $(a_n)\sus M$ with $\lim a_n=x$, we say that {\em $x$ is close to 
$M$}. The theorem says that every uniformly continuous real function defined on a~set of 
fractions $M$ has a~unique limit extension to any real $x$ close to $M$. Our HMC {\em min-max principle} is as follows.

\begin{thm}[HMC min-max principle]\label{thm_minmax}
For every uniformly continuous function $f\cc M\to\R$ defined on a~nonempty bounded set $M\sus\Q$
there exist real numbers $y$ and $y'$ that are close to $M$ and are such that
$$
\forall\,x\in M\,\big(f(y)\le f(x)\le f(y')\big)\;.
$$
\end{thm}
Thus the extended $f$ attains ``on $M$'' at $y$ a~minimum value and at $y'$ a~maximum value. In the displayed formula 
one can replace $M$ with any set $M'$ arising from $M$ by adding to it at most countably many real numbers close to $M$.

The result in \cite{klaz_transc} on which everything hinges is not the min-max principle but a~HMC version of the {\em vanishing derivative 
principle}. One of the standard formulations of it says that if a~function 
$f\cc(a,b)\to\R$ has derivative 
$f'(c)\ne0$, where $a<c<b$ are real numbers,  then $f$ does not have at $c$ local extreme. See \cite{klaz_transc} for our HMC version of this principle .

In CMA real numbers play role of ideal elements which are invoked when they are needed. In 
Theorem~\ref{thm_minmax}, $f$ need not attain extremal value at any element of $M$ but there 
always exist ideal elements $y$ and $y'$, which are Cauchy sequences in $M$, which do the job. Standard 
Cantorean real numbers are equivalence blocks in
$C/\!\!\sim$
where $C$ is the set of rational Cauchy sequences and $\sim$ is the equivalence relation 
given by $(a_n)\sim(b_n)$ iff $\lim\,(a_n-b_n)=0$. We cannot use such real numbers in CMA 
because each of them is uncountable. This cannot be fixed by the axiom of 
choice (AC) by selecting from each equivalence block one representing rational Cauchy 
sequence. Each resulting real number is a~HMC set but AC was applied to an uncountable set of uncountable sets.  We need  real numbers that are HMC from the start.

Such real numbers are well known, they are the Dedekindean real numbers introduced in \cite{dede}. Historically this was the first formalization of real numbers, by means of (Dedekind) {\em cuts} on the set of fractions 
$\Q$. Recall that $X\sus\Q$ is a~{\em cut} if (i) 
$X,\Q\setminus X\ne\emptyset$, (ii) always $a\in\Q$, $b\in X$, $a<b$ $\Rightarrow$ $a\in X$ 
and (iii) $X$ does not have maximum element. But cuts do not capture the required feature 
of real numbers as arbitrarily precise rational approximations, see Theorem~\ref{thm_extend}. 
Therefore in the rest of our article we develop HMC Cantorean real numbers. Also, the arithmetic 
of cuts is a~bit cumbersome. We will proceed in a~quite detailed manner because Cantor's (and 
Heine's and M\'eray's) construction
of real numbers as equivalence blocks of rational Cauchy sequences is well known, but its 
modification that we need in CMA and CNT is, as far as we know, new. 

The belief in indispensability of uncountable sets in Mathematical Analysis is universal. It is 
supported by the fact, often taught in courses of analysis, that the set $\R$ of real 
numbers is uncountable. Typical function in real analysis like $f\cc I\to\R$, where $I\sus\R$ is 
a~nontrivial real interval, is an uncountable set. We regard uncountable sets as 
problematic because almost all of their elements cannot be described by finite means. But 
we also know that for many mathematicians they are their second nature. Individual real numbers, as originally conceived by R.~Dedekind in \cite{dede}, are HMC sets. Also, it is not written 
in stone that in analytical arguments one has to use everything of the mentioned sets $\R$ and $f$, 
maybe some tiny countable parts would suffice for the considered problem. Exactly this we did in \cite{klaz_transc} for the 
transcendence of $\mathrm{e}$. We think that this approach can be extended to many other results in 
Mathematical Analysis and Number Theory, and regard the interest and importance of this undertaking as 
self-evident. 

In Section~2 we briefly review constructions of natural numbers, of the ring of integers and of 
the ordered field of fractions. Section~3 is devoted to the construction of HMC Cantorean 
real numbers and to the proofs that they form a~weak ordered field 
(Theorem~\ref{thm_HMCreals}) and have the weak least upper bound property 
(Theorem~\ref{thm_LUBPforR}). The qualification ``weak'' indicates that in some parts of the 
result the equality relation $=$ is relaxed to the equivalence relation $\sim$. In the last Section~4 we give concluding comments.

\section{Natural numbers, integers, fractions}

We begin with the {\em natural numbers} 
$$
\omega=\{0,\,1,\,2,\,\ds\}
$$ 
where $0=\emptyset$, $1=\{0\}$, 
$2=\{0,1\}$ and so on. More precisely, by the axiom of infinity there exists an inductive set,  
and we define $\omega$ as the intersection of all inductive sets. Then one introduces 
standard addition $+$ and multiplication $\cdot$ on $\omega$ and shows that both (binary) operations 
are commutative and associative, that $\cdot$ is distributive to $+$ and that $0$, resp. $1$, is 
neutral to $+$, resp. $\cdot$\;.

But additive inverses are missing. We set
$$
\Z:=\omega\cup((\omega\setminus\{0\})\times\{0\})
$$
and write, as usual, $-n$ instead of $(n,0)\in\Z$. We set $-0:=0$. We call the elements of 
$\Z$ {\em integers}. One easily extends both operations $+$ and $\cdot$ from $\omega$ to $\Z$. Their previous
properties are preserved and since 
$$
\forall\,n\in\omega\,\big(n+(-n)=0\big)\;,
$$
we get additive inverses. So $(\Z,0,1,+,\cdot)$ is a~commutative ring with identity.

Multiplicative inverses are still missing. We set $Z:=\Z\times(\Z\setminus\{0\})$ 
and write, as usual, $\frac{m}{n}$ or $m/n$ for $(m,n)\in Z$. The {\em identity relation} $\sim$ on $Z$ is
$$
k/l\sim m/n\stackrel{\mathrm{def}}{\iff}kn=lm\;.
$$
It is an equivalence relation on $Z$. Thus we set
$$
\Q:=Z/\!\sim
$$
and call the elements of $\Q$, which are equivalence blocks $[m/n]_{\sim}$, {\em rational 
numbers} or {\em fractions}. Every $\al\in\Q$ is a~countable HMC set and the question if 
$k/l\sim m/n$ is algorithmicly decidable. We will abuse notation and write often, as is common, simply $\frac{m}{n}$ or $m/n$ instead of $[m/n]_{\sim}$. We say that a~fraction $\al\in\Q$ is {\em integral} if $\al=[m/1]_{\sim}$. The map
$$
\Z\ni m\mapsto[m/1]_{\sim}\in\Q
$$
is a~ring isomorphism. 

One easily extends the operations $+$ and $\cdot$ on $\Z$ from integral fractions to $\Q$. All previous properties of $+$ and $\cdot$ are preserved and since 
$$
[m/n]_{\sim}\cdot[n/m]_{\sim}=[mn/nm]_{\sim}=
[1/1]_{\sim}=1_{\Q}\;,
$$
we get multiplicative inverses. Thus $(\Q,0_{\Q},1_{\Q},+,\cdot)$ is 
a~field. It is even an ordered field: if $l,n>0$ (here $>$ is the standard linear order on $\Z$, obtained from the linear order $(\omega,\in)$) then
$$
[k/l]_{\sim}<[m/n]_{\sim}\stackrel{\mathrm{def}}{\iff}kn<lm\;.
$$
One shows that $(\Q,<)$ is a~linear ordering and that 
$$
(\Q,\,0_{\Q},\,1_{\Q},\,+,\,\cdot,\,<)
$$ 
is an ordered field.

One thing is still missing. The ordered field $\Q$ does not have the {\em least upper bound property}. For example, 
the nonempty set
$$
\{\al\in\Q\;|\;\al^2<2\}\sus\Q
$$
is in $<$ bounded from above, but has no least upper bound. 

\section{HMC Cantorean real numbers}

As is well known, in {\em real numbers} the last deficiency is removed. We turn 
to them now. In HMC reals there will be some twists.

Let $X$ be a~set, $+\cc X\times X\to X$ be a~(binary) operation on $X$, $A\sus X\times X$ be 
a~(binary) relation on $X$ and $\sim$ be an equivalence relation on $X$. We say that $+$ is {\em congruent to $\sim$} if for every $a$, $a'$, $b$, $b'$ in $X$ it holds that
$$
(a\sim a'\wedge b\sim b')\Rightarrow a+b\sim a'+b'\;.
$$
Similarly, $A$ is {\em congruent to $\sim$} if for every $a$, $a'$, $b$, $b'$ in $X$ it holds that
$$
(a\sim a'\wedge b\sim b')\Rightarrow(aAb\iff a'Ab')\;.
$$

\begin{defi}[ordered fields congruent to $\sim$]\label{def_weakOF}
Let $X\ne\emptyset$ be a~set and $\sim$ be an equivalence relation on $X$. An ordered field (on $X$) congruent to $\sim$ is a~six-tuple
$$
X_{\mathrm{OF}\sim}:=(X,\,0_X,\,1_X,\,+,\,\cdot,\,<)
$$
such that $(X,\,0_X,\,1_X,\,+,\,\cdot)$ is a~commutative ring with identity, the operations 
$+$ and $\cdot$ on $X$ are congruent to $\sim$, $<$ is an irreflexive and transitive relation on 
$X$ that is congruent to $\sim$, the two ordering axioms hold, namely for every $a,b,c\in X$ one has that
$$
a<b\Rightarrow a+c<b+c\,\text{ and }\,a,\,b>0_X\Rightarrow a\cdot b>0_X\;,
$$
and $X_{\mathrm{OF}\sim}$ has two more properties. First, weak multiplicative inverses exist,
$$
\forall\,a\in X\,\big(a\not\sim 0_X\Rightarrow\exists\,b\in X\,\big(a\cdot b\sim 1_X\big)\big)\;.
$$
Second, $<$ is weakly trichotomic,
$$
\forall\,a,\,b\in X\,\big(a<b\vee b<a\vee a\sim b\big)\;.
$$
\end{defi}
We remind that the requirement on $(X,0_X,1_X,+,\cdot)$ means that $+$ and $\cdot$ are associative and 
commutative, $\cdot$ is distributive to $+$, the element $0_X$ (resp. $1_X$) is neutral to $+$ (resp. 
$\cdot$) and every $a\in X$ has the additive inverse $-a\in X$.

For example, if $=$ is the standard set-theoretic equality, which the axiom of extensionality
characterizes by the equivalence
$$
\text{$x=y$ iff $\forall\,z\,\big(z\in x\iff z\in y\big)$}\;,
$$
then 
$$
\Q_{\mathrm{OF}}=\Q_{\mathrm{OF=}}:=(\Q,\,0/1,\,1/1,\,+,\,\cdot,\,<)
$$
is an ordered field congruent to $=$. This is a~cumbersome way of saying that 
$\Q_{\mathrm{OF}}$ is an ordered field (we defined it in the previous section). Now we define an ordered 
field congruent to an equivalence 
relation weaker than $=$.

Symbols $k$, $l$, $m$, $n$, $n_0$, $n_1$, $\ds$, $n_1'$, $n_2'$, $\ds$ refer to elements of 
$$\N:=\omega\setminus\{0\}\;. 
$$
A~{\em sequence $(a_n)$ in (a~set) $X$} is a~function $a\cc\N\to X$ from
$\N$ to $X$, i.e., a~set of ordered pairs $(a_n)\sus\N\times X$
such that for every $m\in\N$ there is exactly one $y\in X$ with $(m,y)\in(a_n)$. One writes 
$a_m$ for this unique $y$. We denote the set of all sequences in $X$ by $X^{\N}$.

We say that a~sequence $(a_n)$ in $\Q$
is {\em Cauchy} if
$$
\forall\,k\,\exists\,n_0\,\big(m,\,n\ge n_0\Rightarrow|a_m-a_n|\le1/k\big)\;.
$$
We denote the set of all such {\em rational Cauchy sequences} by $C$. The {\em closeness relation} $\sim$ on $C$ is 
$$
(a_n)\sim(b_n)\stackrel{\mathrm{def}}{\iff}
\forall\,k\,\exists\,n_0\,\big(n\ge n_0\Rightarrow|a_n-b_n|\le1/k\big)\;.
$$
Since $(a_n)$ and $(b_n)$ are Cauchy, we can 
equivalently replace the last implication with
$$
m,\,n\ge n_0\Rightarrow|a_m-b_n|\le1/k\;.
$$
It is easy to see that $\sim$ is an equivalence relation on $C$. In the Introduction we mentioned 
that the {\em standard Cantorean real numbers $\R$} are 
$$
\R:=C/\!\sim\;.
$$
They are not HMC sets as every $\al\in\R$ is uncountable. We modify them as follows.

\begin{defi}[HMC reals]\label{def_HMCreals}
We define {\em HMC} real numbers simply by setting
$$
\R:=C\;.
$$
So (our) real numbers are exactly rational 
Cauchy sequences.
\end{defi}
Clearly, every HMC real number is a~HMC set. Their set $C$ is uncountable.

We define arithmetic on $C$ by means of the ordered field $\Q_{\mathrm{OF}}$. Suppose that $(a_n)$ and $(b_n)$ lie in $C$. We set 
$0_C:=(0/1,0/1,\ds)$, $1_C:=(1/1,1/1,\ds)$,
$(a_n)+(b_n):=(a_n+b_n)$, $(a_n)\cdot(b_n):=(a_n\cdot b_n)=(a_nb_n)$ and
$$
(a_n)<(b_n)\stackrel{\mathrm{def}}{\iff}
\exists\,k\,\exists\,n_0\,\big(n\ge n_0\Rightarrow a_n<b_n-1/k\big)\;.
$$
Again, since $(a_n)$ and $(b_n)$ are Cauchy, we can 
equivalently replace the last implication with
$$
m,\,n\ge n_0\Rightarrow a_m<b_n-1/k\;.
$$
The notation $(a_n)\lesssim(b_n)$ means that $(a_n)<(b_n)$ or $(a_n)\sim(b_n)$, and similarly for $\gtrsim$. We show that $(C,0_C,1_C,+,\cdot,<)$ is an ordered field congruent to $\sim$. Its ring 
structure is immediate from the following more general construction.

\begin{prop}[$\N$-th powers of rings]\label{prop_powerRing}
$(R,0_R,1_R,+,\cdot)$ is a~commutative ring with identity and $P:=R^{\N}$. Then 
$$
P_{\mathrm{R}}:=(P,\,0_P,\,1_P,\,+,\,\cdot)\;,
$$
where $0_P:=(0_R,0_R,\ds)$, $1_P:=(1_R,1_R,\ds)$ and the operations $+$ and $\cdot$ on $P$ are defined component-wisely from those on $R$, is a~commutative ring with identity.
\end{prop}
\proof
Satisfaction of the axioms of a~commutative ring with identity in $P_{\mathrm{R}}$ is immediate because they hold in every component. 
\eproof

CMA views HMC reals as follows.

\begin{thm}[HMC reals form a~weak ordered field]\label{thm_HMCreals}
The structure
$$
\R:=(C,\,0_C,\,1_C,\,+,\,\cdot,\,<)
$$
defined above is an ordered field congruent to the closeness relation $\sim$, in the sense of Definition~\ref{def_weakOF}.
\end{thm}
\proof
It is clear that $0_C$ and $1_C$ are in $C$. Let $(a_n)$ and $(b_n)$ lie in $C$. Clearly, $(a_n)+(b_n)=(a_n+b_n)\in C$. We treat $(a_n)\cdot (b_n)=(a_n b_n)$ in more detail. For a~given $k$ there 
is an $n_0$ such that $m,n\ge n_0$ $\Rightarrow$ $|a_m-a_n|,|b_m-b_n|\le1/k$. Hence there is an $l$
(independent of $k$) such that $\forall\,n\,\big(|a_n|,|b_n|\le l\big)$.
Thus for every $m,n\ge n_0$,
$$
|a_mb_m-a_nb_n|\le|a_m|\cdot|b_m-b_n|+|a_m-a_n|\cdot|b_n|\le 2l/k
$$
and we see that $(a_n)\cdot (b_n)\in C$. One can similarly prove that $+$ and $\cdot$ are 
congruent to $\sim$. We have shown that $C\sus\Q^{\N}$ contains $0_C$ and $1_C$ and 
is closed to the operations $+$ and $\cdot$. By Proposition~\ref{prop_powerRing}, 
$(C,0_C,1_C,+,\cdot)$ is a~commutative ring with identity. We show that it has weak multiplicative inverses. If $(a_n)\in C$ with $(a_n)\not\sim 0_C$ then $|a_n|\ge1/k$ for every $n\ge n_0$ and some $k$. We define $(b_n)\in\Q^{\N}$ by
$$
b_n:=\left\{
\begin{array}{lll}
0 & \ds & a_n=0\,\text{ and} \\
1/a_n & \ds & a_n\ne0\;. 
\end{array}
\right.
$$
Since for a~given $l$ there is an $n_1$ such that $m,n\ge n_1$ $\Rightarrow$ $|a_m-a_n|\le1/l$ and we can also assume that $n\ge n_1$ $\Rightarrow$ $|a_n|\ge1/k$, for every $m,n\ge n_1$ it holds that
$$
|b_m-b_n|=\bigg|\frac{1}{a_m}-\frac{1}{a_n}\bigg|=
\frac{|a_n-a_m|}{|a_m|\cdot|a_n|}\le\frac{k^2}{l}
$$
and $(b_n)\in C$. Since $a_n b_n=1/1$ for every $n\ge n_0$, we see that 
$(a_n)\cdot(b_n)\sim 1_C$.

We verify the properties of $\R$ concerning $<$.  Clearly, $<$ is irreflexive. If $(a_n)<(b_n)$ and 
$(b_n)<(c_n)$ then there exist $k$ and $n_0$ such that for every $n\ge n_0$,
$$
a_n<b_n-1/k\,\text{ and }\,b_n<c_n-1/k\;.
$$
Hence $a_n<c_n-2/k<c_n-1/k$ for every $n\ge n_0$ and $<$ is transitive. We show that $<$ is 
congruent to $\sim$. Suppose that $(a_n)$, $(a_n')$, $(b_n)$ and $(b_n')$ lie in $C$, $(a_n)<(b_n)$,
$(a_n)\sim(a_n')$ and $(b_n)\sim(b_n')$. Then there exist $k$ and $n_0$ such that 
$$
n\ge n_0\Rightarrow a_n<b_n-1/k\;.
$$
Since $(a_n)\sim(a_n')$ and $(b_n)\sim(b_n')$, there exist an $n_1\ge n_0$ such that
$$
n\ge n_1\Rightarrow a_n'<b_n'-1/2k\;.
$$
Hence $(a_n')<(b_n')$.

Suppose that 
$(a_n),(b_n)\in C$ with $(a_n)\not\sim(b_n)$. Then there is a~$k$ such that $|a_n-b_n|>1/k$ for 
infinitely many $n$. Thus $a_n<b_n-1/k$ for 
infinitely many $n$ or $b_n<a_n-1/k$ for 
infinitely many $n$. Since $(a_n),(b_n)\in C$, in the former case there is an $n_0$ such that 
$n\ge n_0$ $\Rightarrow$ $a_n<b_n-1/2k$ and 
$(a_n)<(b_n)$. In the latter case the same argument gives that $(b_n)<(a_n)$. We have shown 
that $<$ is weakly trichotomic. Let $(a_n)$, $(b_n)$ and $(c_n)$ lie in $C$. If $(a_n)<(b_n)$ then there exist $k$ and $n_0$ such that
$$
n\ge n_0\Rightarrow a_n<b_n-1/k\;.
$$
Thus $a_n+c_n<b_n+c_n-1/k$ for every $n\ge n_0$ and $(a_n)+(c_n)<(b_n)+(c_n)$. Similarly, if $(a_n),(b_n)>0_C$ then there exist $k$ and $n_0$ such that
$$
n\ge n_0\Rightarrow 1/k<a_n,\,b_n\;.
$$
Thus $1/k^2<a_n b_n$ for every $n\ge n_0$ and $(a_n)\cdot(b_n)>0_C$. This proves the two 
ordering axioms for $\R$ and concludes the proof of the theorem.
\eproof

\noindent
Before we turn to the proof of the least upper bound property for $\R$ we have to clarify how $\Q$
is contained in $\R$. The situation is actually similar to the containment of $\Z$ 
in $\Q$. We call the constant sequences 
$${\textstyle
E(\frac{m}{n}):=(\frac{m}{n},\,\frac{m}{n},\,\ds)\in C,\ m/n\in\Q\;, 
}
$$
{\em 
rational\, {\em HMC} reals} and denote their (countable) set by $\Q_C$. The structure of the ordered field $\R$ 
congruent to $\sim$ restricts on $\Q_C$ to the structure of an ordinary ordered field (i.e., on 
$\Q_C$ the relation $\sim$ upgrades to $=$). The map $E\cc\Q\to\Q_C$ is then an isomorphism of ordered fields.

We show that, unlike $\Q_{\mathrm{OF}}$, HMC reals have the least upper bound property. In the weak sense, though, with $=$ relaxed to $\sim$. In CMA one can use only at 
most countable subsets of $C$, but the result holds for any subset and we prove it as such.

\begin{thm}[$\R$ has weak LUBP]\label{thm_LUBPforR}
{\em HMC} real numbers have the weak least upper bound property. Namely, for every nonempty set $B\sus C$ if $(b_n)\lesssim(a_n)$ for every $(b_n)\in B$ and some $(a_n)\in C$, then $B$ has a~least upper bound. It is a~($\sim$-unique) sequence $(a_n')\in C$ such that
\begin{itemize}
    \item $(b_n)\lesssim(a_n')$ for every $(b_n)\in B$ and
    \item for every $(c_n)\in C$ with $(c_n)<(a_n')$ there is a~$(b_n)\in B$ with $(c_n)<(b_n)$. 
\end{itemize}
\end{thm}
\proof
Suppose that $B\sus C$ is a~nonempty set and $(a_n)\in C$ is an upper bound of $B$. Clearly, we may take $(a_n)$ to be 
$E(m/1)$ for some $m\in\N$. In other words, $\R$ is Archimedean. In the following procedure with four commands we inductively define two rational sequences $(a_n')$ and $(b_n)$ in $\Q$.
\begin{enumerate}
    \item (initialization) $a_1':=m/1$ and $b_1:=1/1$.
    \item (branching) Suppose that the fractions $a_1'$, $\ds$, $a_n'$ and $b_1$, $\ds$, $b_n$ have been defined. Is $E(a_n'-b_n)$ still an upper bound of $B$?
    \item If YES, set $a_{n+1}':=a_n'-b_n$, $b_{n+1}:=b_n$ and go back to command $2$.
    \item If NO, set $a_{n+1}':=a_n'$, $b_{n+1}:=\frac{1}{1+1/b_n}$ and go back to command $2$.
\end{enumerate}
The sequence $(a_n')$ is non-increasing and for every $n$, $E(a_n')$ is an upper bound
of $B$. We show that $(a_n')\in C$ and is the desired least upper bound of $B$.

Clearly, command~$4$ is performed infinitely many times. Thus the sequence
$${\textstyle
(b_n)=\big(\frac{1}{1},\,\frac{1}{1},\,\ds,\,\frac{1}{1},\,\frac{1}{2},\,\frac{1}{2},\,\ds,\,\frac{1}{2},\,\frac{1}{3},\,\frac{1}{3},\,\ds,\,\frac{1}{3},\,\frac{1}{4},\,\ds,\,\ds,\,\ds\big)
}
$$
and goes to $\frac{0}{1}$. We denote by $1\le m_1<m_2<\ds$ those steps $n=m_i$ in the procedure 
when command~$4$ is performed, and select elements 
$$
b_i'=(d_n^i)\in B
$$ 
such that in step $n=m_i$ one 
has that $E(a_n'-b_n)<b_i'$. The last inequality means that there exist $l_i$ and $n_i'$ in $\N$ such that $n\ge n_i'$ $\Rightarrow$ 
$a_{m_i}-1/i<d_n^i-1/l_i$. Then for every $i\in\N$,
$$
n\ge m_i\Rightarrow E(a_{m_i}')\ge E(a_n')\gtrsim b_i'>E(a_{m_i}'-1/i)\;.
$$
Thus for every $i$ for large $n$ one has that $a_{m_i}'\ge a_n'>a_{m_i}'-2/i$ and $(a_n')\in C$.

We show that $(a_n')$ is an upper bound of $B$. Suppose for the contradiction that 
$(a_n')<b$ for some $b=(d_n)\in B$. Then there exist $k$ and $n_0$ such that
$$
n\ge n_0\Rightarrow a_n'<d_n-1/k\;.
$$
Since $(a_n')\in C$, there is an $n_1$ such that
$m,n\ge n_1$ $\Rightarrow$ $|a_m'-a_n'|\le 1/2k$. But then 
$$
n\ge N:=\max(\{n_0,\,n_1\})\Rightarrow
a_N'<d_n-1/2k\;.
$$
Thus $E(a_N')<b$, a~contradiction. It remains to show that $(a_n')$ is the least upper bound 
of $B$. Let $(c_n)\in C$ be any sequence with $(c_n)<(a_n')$. Thus there exist $k$ and $n_0$ such that 
$$
n\ge n_0\Rightarrow c_n<a_n'-1/k\;.
$$
But then for every $n\ge\max(\{n_0,m_k\})$, 
$$
E(c_n)<E(a_n'-1/k)\le E(a_{m_k}'-1/k)<b_k'\;.
$$
Thus for every $n\ge\max(\{n_0,m_k,n_k'\})$ one has that 
$$
c_n<d^k_n-1/l_k\;.
$$
So $(c_n)<b_k'$ and $(c_n)$ is not an upper bound of $B$.
\eproof

\section{Concluding remarks}\label{sec_conclRema}

$\R$ is not just an ordinary set in set theory, it is a~set that in a~sense gave birth to set theory, 
and rightly \cite{jech_ST} devotes Chapter~4 and ten pages to it. There are several 
modern books on real numbers, of which we explicitly mention only \cite{buko} and 
\cite{stil}. In this article we lighted a~facet of real numbers that is not considered in 
these books. 

In \cite{klaz_transc} we developed a~fragment of Countable Mathematical Analysis and Countable 
Number Theory. The present 
version~3 of \cite{klaz_transc} has to be revised and take into account the treatment of real numbers in this article. We apologize to 
the readers of \cite{klaz_transc} for this deficiency and hope to produce the corresponding 
revision soon.

\medskip\noindent
{\em Department of Applied Mathematics\\
Faculty of Mathematics and Physics\\
Charles University\\
Malostransk\'e n\'am\v est\'\i\ 25\\
118 00 Praha\\
Czechia\\
{\tt klazar@kam.mff.cuni.cz}
}


\begin{thebibliography}{20}

\bibitem{buko}
L. Bukovsk\'y, {\em The Structure of the Real Line}, Springer, Basel 2011

\bibitem{dede}
R. Dedekind, {\em Was sind und was sollen die Zahlen?}, Vieweg und Sohn, Braunschweig 1872.

\bibitem{jech_ST}
T. Jech, {\em Set Theory}, Springer-Verlag, Berlin 2003

\bibitem{klaz_transc}
M. Klazar, A chapter in Countable Number Theory: the transcendence of Euler's number, 	arXiv:2301.08142 [math.LO], 2023, 42 pp. 

\bibitem{stil}
J. Stillwell, {\em The Real Numbers. An Introduction to Set Theory and Analysis}, Springer, Cham 2013

\end{thebibliography}
\end{document}